\documentclass[a4paper, 11pt, oneside, onecolumn,sort&compress]{elsarticle}
\usepackage{amssymb}
\usepackage{mathrsfs}
\usepackage{amsfonts}
\usepackage{ntheorem}

\newtheorem*{thmA}{Theorem A}
\newtheorem*{thmB}{Theorem B}
\newtheorem*{thmC}{Theorem C}
\newtheorem{theorem}{Theorem}[section]
\newtheorem{lemma}{Lemma}[section]

\newdefinition{example}{Example}[section]
\newdefinition{remark}{Remark}[section]
\newdefinition{definition}{Definition}[section]
\newproof{proof}{Proof}[section]
\newproof{pot}{Proof of Theorem \ref{the1.1}}
\usepackage{amssymb,amsmath}

\numberwithin{equation}{section} \journal{Bulletin des Sciences Math\'{e}matiques }

\begin{document}
\begin{frontmatter}

\title{\textbf{Multiple solutions for Schr\"{o}dinger equations under the partially sublinear case}
\tnoteref{label1}} \tnotetext[label1]{This work was supported by the
National Natural Science Foundation of China (Grant Nos. 11271299
and 11001221), the Mathematical Tianyuan Foundation of China (No.
11126027), Natural Science Foundation of Shanxi Province
(2014021009-1) and the Shanxi University Yinjin Rencai Jianshe
Xiangmu (No. 010951801004).}
\author  {Xiaojing Feng\corref{cor1}}
\ead{fengxj@sxu.edu.cn}

\cortext[cor1]{Corresponding author.}

\address{School of Mathematical Sciences, Shanxi University, Taiyuan 030006, People's Republic of China}

\begin{abstract}

In this paper, we prove the infinitely many solutions to a class of sublinear
Schr\"{o}dinger-Poisson equations by using an extension of Clark's theorem
established by Zhaoli Liu and Zhi-Qiang Wang.
\end{abstract}
\end{frontmatter}

\section{Introduction and main results}
In this paper we study the existence and
multiplicity of solutions for the following Schr\"{o}dinger-Poisson system:
$$\left\{\begin{array}{ll}
-\Delta u+u+\phi u=K(x)f(x,u),&{\rm in}\ \mathbb{R}^3,\\-\Delta \phi=u^2,\
\lim_{|x|\to \infty}\phi(x)=0,&{\rm in}\ \mathbb{R}^3.
\end{array}\right.\eqno{(1.1)}$$

There have been a lot of research for system (1.1), after it was first introduced in \cite{bf} as a physical
model describing a charged wave interacting with its own electrostatic field in quantum mechanic.
Recently, many authors were interested in the case $f(x,u)=|u|^{p-1}u,\ 1 < p < 5$(see\cite{c,d}).
In \cite{c}, the author proves the existence of a radial solution when
$3 < p < 5$. On the other hand, the non-radial solution is considered in \cite{d}.
J. Seok \cite{s} proved
that (1.1) had infinitely many high energy solutions under the condition that $f(x,u)$
is superlinear at infinity in $u$.

In \cite{sun}, Sun studied the existence of
solutions of the problem

$$\left\{\begin{array}{ll}
-\Delta u+V(x)u+\phi u=f(x,u),&{\rm in}\ \mathbb{R}^3,\\-\Delta \phi=u^2,
\ \lim_{|x|\to \infty}\phi(x)=0,&{\rm in}\ \mathbb{R}^3.
\end{array}\right.\eqno{(1.2)}$$

Under the case where
$F(x,u)$ is a special case of subquadratic at infinity and obtained the following theorem.

\begin{thmA}\cite{sun} Assume that the following conditions hold:

($A_1$) $V\in C(\mathbb{R}^3,\mathbb{R})$ and $\inf_{x\in \mathbb{R}^3}V(x)\geq\beta>0$;

($A_2$) For every $M>0$, $meas\{x\in \mathbb{R}^3:\ V(x)\leq M\}<\infty$, where meas$(\cdot)$ denotes
the Lebesgue measure in $\mathbb{R}^3$;

($A_3$) $F(x,u)=b(x)|u|^{p+1}$, where $F(x,u)=\int_0^uf(x,y)dy,\ b:\mathbb{R}^3\to \mathbb{R}^+$
is a positive continuous function such that $b\in L^{2/(1-p)}(\mathbb{R}^3)$ and $0<p<1$ is a constant.\\
Then (1.2) possesses infinitely many solutions.
\end{thmA}
In this Theorem, as $F(x,u)$ is special case, many sublinear functions in mathematical physics do not
satisfy the condition ($A_3$). Soon after, the authors use the genus properties in critical theory to
generalize Theorem A by relaxing assumption ($A_3$) and obtained the following theorem.

\begin{thmB}\cite{ct} Assume that $V$ and $f$ satisfy $A_1$, $A_2$ and the following conditions:

(B$_1$) $f\in C(\mathbb{R}^3\times \mathbb{R},\mathbb{R})$ and
and there exist constants $1<\gamma_1<\gamma_2<\ldots<\gamma<2$ and functions
$a_i\in L^{2/(2-\gamma_i)}(\mathbb{R}^3,[0,+\infty))(i=1,2,\ldots,m)$ such that
$$|f(x,z)|\leq\sum_{i=1}^ma_i(x)|z|^{\gamma_i-1},\ \forall\ (x,z)\in\mathbb{R}^3\times \mathbb{R};$$

($B_2$) There exists an open set $\Lambda\subset\mathbb{R}^3$ and three constants $\delta>0$, $\gamma_0\in(1,2)$
and $\eta>0$ such that
$$F(x,z)\geq\eta|z|^{\gamma_0},\ \forall\ (x,z)\in \Lambda\times[-\delta,\delta];$$

($B_3$) $f(x,-z)=-f(x,z),\ (x,z)\in \mathbb{R}^3\times \mathbb{R}.$.\\
Then (1.2) possesses infinitely many solutions.
\end{thmB}

Recently, the authors obtained an extension of Clark's theorem as follows.

\begin{thmC}\cite{lw}
Let $X$ be a Banach space, $\Phi\in C^1(X,\mathbb{R})$. Assume $\Phi$ is even and satisfies the (PS) condition,
bounded from below, and $\Phi(0)=0$. If for any $k\in \mathbb{N}$, there exists a $k-$dimensional subspace $X^k$
of $X$ and $\rho_k>0$ such that $\sup_{X^k\cap S_{\rho_k}}\Phi<0$, where $S_\rho=\{u\in X|\|u\|=\rho\}$, then
at least one of the following conclusions holds.

(i) There exists a sequence of critical points $\{u_k\}$ satisfying $\Phi(u_k)<0$ for all $k$ and
$\|u_k\|\to 0$ as $k\to \infty$.

(ii) There exists $r>0$ such that for any $0<a<r$ there exists a critical point $u$ such that $\|u\|=a$
and $\Phi(u)=0$.
\end{thmC}

In this paper, we will generalize Theorema A and B by using the Theorem C. Our main result is as follows.
\begin{theorem}\label{the1.1}
Assume that $f$ satisfies ($B_3$) and the following conditions:

($f_1$) There exist $\delta>0,\ 1\leq \gamma<2,\ C>0$ such that $f\in
C(\mathbb{R}^3\times[-\delta,\delta],\mathbb{R})$ and
$|f(x,z)|\leq C|z|^{\gamma-1}$;

($f_2$) $\lim_{z\to 0}F(x,z)/|z|^2=+\infty$
uniformly in some ball $B_r(x_0)\subset\mathbb{R}^3$, where $F(x,z)=\int_0^zf(x,s)ds$.

($f_3$) $K:\mathbb{R}^3\to \mathbb{R}^+$
is a positive continuous function such that $K\in L^{2/(1-\gamma)}(\mathbb{R}^3)\cap L^{\infty}(\mathbb{R}^3)$.

Then (1.1) possesses infinitely many solutions $\{u_k\}$ such that
$\|u_k\|_{L^\infty}\to 0$ as $k\to \infty$.
\end{theorem}
\begin{remark}
It is to see that ($B_2$) implies ($f_2$) and ($f_1$) is weaker than ($B_1$) in essence.
\end{remark}

\begin{remark} Throughout the paper we denote by $C>0$ various positive constants which may
vary from line to line and are not essential to the problem.
\end{remark}
The paper is organized as follows: in Section 2, some preliminary results are presented.
Section 3 is devoted to the proof of Theorem 1.1.

\section{Preliminary\label{Section 2}}

In this Section, we will give some notations and Lemma that will be used throughout this paper.

Let $H^1=H^1(\mathbb{R}^3)$ be the completion of $C_0^\infty(\mathbb{R}^3)$ with respect to the norm
$$\|u\|^2=\int_{\mathbb{R}^3}|\nabla u|^2+u^2dx.$$
Moreover, we denote the completion of $C_0^\infty(\mathbb{R}^3)$ with respect to the norm
$$\|u\|_{D^1}=\int_{\mathbb{R}^3}|\nabla u|^2dx$$
by $D^1=D^1(\mathbb{R}^3)$.
To avoid lack of compactness, we need consider the set of radial functions as follows:
$$H=H_r^1(\mathbb{R}^3)=\{u\in H^1(\mathbb{R}^3)| u(x)=u(|x|)\}$$
and
$$D=D_r^1(\mathbb{R}^3)=\{u\in D^1(\mathbb{R}^3)| u(x)=u(|x|)\}.$$
Here we note that the continuous embedding $H_r^1(\mathbb{R}^3)\hookrightarrow L^q(\mathbb{R}^3)$
is compact for any $q\in (2,6)$.

We observe that by the Lax-Milgram theorem, for given $u\in H^1$, there exists a unique solution
$\phi=\phi_u\in D^1$ satisfying $-\Delta \phi_u=u^2$ in a weak sense. The function $\phi_u$ is represented by
$$\phi_u(x)=\frac{1}{4\pi}\int_{\mathbb{R}^3}\frac{u^2(y)}{|x-y|}dy,$$
and it has the following properties.

\begin{lemma}[\cite{r}] The following properties hold:

(i) there exists $C>0$ such that for any $u\in H^1(\mathbb{R}^3)$,
$$\|\phi_u\|_{D^1}\leq C\|u\|^2_{L^{12/5}},\ \int_{\mathbb{R}^3}|\nabla \phi_u|^2dx
=\int_{\mathbb{R}^3}\phi_u u^2dx\leq C\|u\|^4;$$

(ii) $\phi_u\geq 0$ for all $u\in H^1$;

(iii) if $u$ is radially symmetric, then $\phi_u$ is radial;

(iv) $\phi_{tu}=t^2\phi_u$ for all $t>0$ and $u\in H^1$;

(v) if $u_j\rightharpoonup u$ weakly in $H$, then, up to a subsequence, $\phi_{u_j}\to \phi_u$ in $D$ and
$$\int_{\mathbb{R}^3}\phi_{u_j}u_j^2dx\to\int_{\mathbb{R}^3}\phi_{u}u^2dx.$$
\end{lemma}

\section{Proofs of the main result \label{Section 3}}

{\bf Proof of Theorem 1.1.} Choose $\hat{f}\in C(\mathbb{R}^N\times\mathbb{R}, \mathbb{R})$ such that $\hat{f}$ is
odd in $u\in \mathbb{R}$, $\hat{f}(x,u)=f(x,u)$ for $x\in \mathbb{R}^N$ and $|u|<\delta/2$, and
$\hat{f}(x,u)=0$ for $x\in \mathbb{R}^N$ and $|u|>\delta$. In order to obtain solutions of (1.1) we consider

$$\left\{\begin{array}{ll}
-\Delta u+u+\phi u=K(x)\hat{f}(x,u),&{\rm in}\ \mathbb{R}^3,\\-\Delta
\phi=u^2,\ \lim_{|x|\to \infty}\phi(x)=0,&{\rm in}\ \mathbb{R}^3.
\end{array}\right.\eqno{(3.1)}$$

Moreover, (3.1) is variational and its solutions are the critical points of the functional defined in
$H$ by
$$J(u)=\frac{1}{2}\|u\|^2+\frac{1}{4}\int_{\mathbb{R}^3}\phi_u(x)u^2dx-\int_{\mathbb{R}^3}K(x)\hat{F}(x,u)dx.$$
From ($f_1$), it is easy to check that $J$ is well defined on $H$ and
$J\in C^1(H^1(\mathbb{R}^3),\mathbb{R})$ (see \cite{ct} for more detail), and
$$J'(u)v=\int_{\mathbb{R}^3}\nabla u\nabla vdx+\int_{\mathbb{R}^3}\phi_u(x)uvdx
-\int_{\mathbb{R}^3}K(x)\hat{f}(x,u)vdx,\ v\in H.$$

Note that $J$ is even, and $J(0)=0$.
For $u\in H^1(\mathbb{R}^3)$,
$$\int_{\mathbb{R}^3}K(x)|\hat{F}(x,u)|dx\leq C\int_{\mathbb{R}^3}K(x)|u|^\gamma dx
\leq
C\|K\|_{L^{\frac{2}{2-\gamma}}(\mathbb{R}^3)}\|u\|^\gamma_{L^{2}(\mathbb{R}^3)}\leq
C\|u\|^\gamma.$$ Hence, it follows from Lemma 2.1 that
$$J(u)\geq\frac{1}{2}\|u\|^2-C\|u\|^\gamma,\ u\in H.\eqno(3.2)$$
We now use the same ideas to prove the (PS) condition. Let $\{u_n\}$ be a sequence in $H$ so that
$J(u_n)$ is bounded and $J'(u_n)\to 0$. We shall prove that $\{u_n\}$ converges. By (3.2), we claim that $\{u_n\}$
is bounded.
Assume without loss of generality that $\{u_n\}$ converges to $u$ weakly in $H$. Observe that
\begin{equation*}
\begin{split}
\|u_n-u\|^2&=\langle J'(u_n)-J'(u),u_n-u\rangle-\int_{\mathbb{R}^3}(\phi_{u_n}u_n-\phi_{u}u)(u_n-u)\\
&-\int_{\mathbb{R}^3}K(x)(\hat{f}(x,u_n)-\hat{f}(x,u))(u_n-u)dx\\
&\equiv I_1+I_2+I_3,
\end{split}
\end{equation*}
It is clear that $I_1\to 0$ as $n\to \infty.$ In the following, we will estimate $I_2$, by computation,
\begin{equation*}
\begin{split}
&|\int_{\mathbb{R}^3}(\phi_{u_n}u_n-\phi_{u}u)(u_n-u)dx|\\
&\leq|\int_{\mathbb{R}^3}\phi_{u_n}(u_n-u)(u_n-u)dx|+|\int_{\mathbb{R}^3}(\phi_{u_n}-\phi_{u})u(u_n-u)dx|\\
&\leq \|\phi_{u_n}\|_{L^6}\|u_n-u\|_{L^{12/5}}^2\\
&+\|\phi_{u_n}-\phi_u\|_{L^6}\|u\|_{L^{12/5}}\|u_n-u\|_{L^{12/5}}\to 0
\end{split}
\end{equation*}
because $\phi_{u_n}\to \phi_u$ in $L^6$ and $u_n\to u$ in $L^{12/5}$ up to a subsequence.

We estimate $I_3$, by using (f$_3$), for any $R>0$,
\begin{equation*}
\begin{split}
&\int_{\mathbb{R}^3}K(x)|\hat{f}(x,u_n)-\hat{f}(x,u)||u_n-u|dx\\
&\leq C\int_{\mathbb{R}^3\setminus B_R(0)}K(x)(|u_n|^\gamma+|u|^\gamma)dx+C\int_{B_R(0)}(|u_n|^{\gamma-1}+|u|^{\gamma-1})|u_n-u|dx\\
&\leq C\left(\|u_n\|^\gamma_{L^2(\mathbb{R}^3\setminus B_R(0))}+\|u\|^\gamma_{L^2(\mathbb{R}^3\setminus B_R(0))}\right)
\|K\|_{L^{\frac{2}{2-\gamma}}(\mathbb{R}^3\setminus B_R(0))}\\
&+C\left(\|u_n\|^{\gamma-1}_{L^{\gamma}(B_R(0))}+\|u\|^{\gamma-1}_{L^{\gamma}(B_R(0))}\right)\|u_n-u\|_{L^{\gamma}(B_R(0))}\\
&\leq C\|a\|_{L^{\frac{2}{2-\gamma}}(\mathbb{R}^3\setminus B_R(0))}+C\|u_n-u\|_{L^{\gamma}(B_R(0))},
\end{split}
\end{equation*}
which implies
$$\lim_{n\to+\infty}\int_{\mathbb{R}^3}K(x)|\hat{f}(x,u_n)-\hat{f}(x,u))||u_n-u|dx=0.$$
Therefore, $\{u_n\}$ converges strongly in $H$ and the (PS) condition holds for $J$. By (f$_2$) and (f$_3$),
for any $L>0$, there exists
$\delta=\delta(L)>0$ such that if $u\in C_0^\infty(B_r(x_0))$ and $|u|_\infty<\delta$
then $K(x)\hat{F}(x,u(x))\geq L|u(x)|^2$,
and it follows from Lemma 2.1 that
$$J(u)\leq\frac{1}{2}\|u\|^2+\frac{1}{4}\|u\|^4-L\|u\|^2_{L^2(\mathbb{R}^3)}.$$
This implies, for any $k\in \mathbb{N}$, if $X^k$ is a $k-$dimensional subspace of $C_0^\infty(B_r(x_0))$
and $\rho_k$ is sufficiently small then $\sup_{X^k\cap S_{\rho_k}}J(u)<0$, where $S_\rho=\{u\in \mathbb{R}^3| \|u\|=\rho\}$.
Now we apply Theorem C to obtain infinitely many solutions $\{u_k\}$ for (3.1) such that

$$\|u_k\|\to 0,\ k\to \infty.\eqno(3.3)$$

Finally we show that $\|u_k\|_{L^\infty}\to 0$ as $k\to \infty$. Let $u$ be a solution of (3.1) and
$\alpha>0$. Let $M>0$ and set $u^M(x)=max\{-M,min\{u(x),M\}\}$. Multiplying both sides of (3.1) with
$|u^M|^\alpha u^M$ implies
$$\frac{4}{(\alpha+2)^2}\int_{\mathbb{R}^3}|\nabla|u^M|^{\frac{\alpha}{2}+1}|^2dx
\leq C\int_{\mathbb{R}^3}|u^M|^{\alpha+1}dx.$$
By using the iterating method in \cite{lw}, we can get the following estimate
$$\|u\|_{L^\infty(\mathbb{R}^3)}\leq C_1\|u\|^{\nu}_{L^6(\mathbb{R}^3)},$$
where $\nu$ is a number in $(0,1)$ and $C_1>0$ is independent of $u$ and $\alpha$. By (3.3) and
Sobolev imbedding Theorem\cite{w},
we derive that $\|u_k\|_{L^\infty(\mathbb{R}^3)}\to 0$ as $k\to \infty$. Therefore, $u_k$ are the solutions of (1.1) as
$k$ sufficiently large. The proof is completed.
\hfill$\Box$

\section{Acknowledgement \label{Section 4}}

 Authors would like to express their sincere gratitude to one anonymous referee for his/her
constructive comments for improving the quality of this paper.

\end{document}